\documentclass[12pt]{amsart}
\usepackage{amscd, amssymb,latexsym,amsmath, amscd, amsmath}
\usepackage[all]{xy}

\newtheorem{theorem}{Theorem}[section]

\newtheorem{prop}[theorem]{Proposition}
\newtheorem{remark}{Remark}

\newtheorem{cor}[theorem]{Corollary}
\newtheorem{conj}{Conjecture}

\begin{document}

\title[Distinguished representations]
{Distinguished non-Archimedean representations}
\author{U. K. Anandavardhanan}

\address{Tata Institute of Fundamental 
Research, Homi Bhabha Road, Bombay - 400 005, INDIA.}
\email{anand@math.tifr.res.in}

  
\date{}

\begin{abstract}
For a symmetric space $(G,H)$, one is interested in understanding the 
vector space of $H$-invariant linear forms on a representation $\pi$ of $G$.
In particular an important question is whether or not the dimension of this 
space is bounded by one. We cover the known results for the pair
$(G={\rm R}_{E/F}{\rm GL}(n),H={\rm GL}(n))$, and then discuss 
the corresponding ${\rm SL}(n)$ case. In this paper, we show that 
$(G={\rm R}_{E/F}{\rm SL}(n),H={\rm SL}(n))$ is a Gelfand pair when
$n$ is odd. When $n$ is even, the space of $H$-invariant forms on $\pi$
can have dimension more than one even when $\pi$ is supercuspidal.
The latter work is joint with Dipendra Prasad.
\end{abstract}

\maketitle

\section{Introduction}

Let $G$ be a group, and $H$ the group of fixed points of an involution on $G$.
A representation $\pi$ of $G$ is said to be distinguished with respect to $H$, 
if the space of $H$-invariant linear forms on $\pi$ is nonzero. More 
generally, if the space ${\rm Hom}_H(\pi,\chi)$ is nonzero for a character 
$\chi$ of $H$, we say that $\pi$ is $\chi$-distinguished with respect to $H$.
In this paper we are interested in the case when $(G,H)$ is defined over a 
$p$-adic field. 

The initial impetus for much of the research in this field came from the work
of Harder, Langlands, and Rapoport \cite{hlr} where they introduce 
the notion of distinguishedness (in terms of the non-vanishing of a certain 
period integral) when $(G,H)$ is defined over the adeles of a number field.
Specifically a careful analysis of distinguishedness for the pair 
$G={\rm R}_{F/{\Bbb Q}}{\rm GL}(2),H={\rm GL}(2)$, where $F$ 
is a real quadratic 
extension of ${\Bbb Q}$, was required to settle the Tate conjecture for
the Hilbert modular surface associated to $G$. 

The philosophy (due to Jacquet) is that distinguishedness for the pair $(G,H)$
often characterise the image of a lift from a suitable group $H^\prime$. For
instance if $G={\rm GL}_n(E)$ and $H={\rm GL}_n(F)$, where $E/F$ a 
quadratic extension
of $p$-adic (or number) fields, the $H$-distinguished representations
are the ones in the image of the base change map from a suitable unitary
group \cite{flicker1}. Conversely if $H$ is a unitary group
(with respect to $E/F$), then $H^\prime$ is ${\rm GL}_n(F)$ 
(see the series of papers
starting with \cite{jacquet1} and \cite{jacquet2}).

From now on let $E/F$ be a quadratic extension of $p$-adic fields.
The groups that we consider will be defined over $F$.
For a representation $\pi$ of $G$, one is interested in a better 
understanding of the space ${\rm Hom}_H(\pi,1)$. A symmetric pair $(G,H)$
is said to have the multiplicity one property (equivalently $(G,H)$ is called 
a Gelfand pair) if the above space has dimension less than or equal to one
for all irreducible admissible representations $\pi$ of $G$.
The pair $({\rm GL}_n(E),{\rm GL}_n(F))$ is known to have the multiplicity 
one property (\cite{flicker1},\cite{hakim1}). Here it is not very hard to 
deduce the multiplicity one property, as it is relatively easy to show
that all the double cosets of $H$ in $G$ are fixed by an involution. 
Nevertheless proving that a pair $(G,H)$ is a Gelfand pair can turn out to be 
quite hard. For instance this is the case when $G={\rm GL}_{2n}(F)$ and 
$H={\rm GL}_n(F)\times {\rm GL}_n(F)$ \cite{jacquet3}. 

An example of a non-Gelfand pair is obtained by taking $G={\rm GL}_n(E)$ and 
$H={\rm U}(n,E/F)$, the quasi-split unitary group with respect to $E/F$. 
But in this
case, $\dim_{\Bbb C}{\rm Hom}_H(\pi,1)$ is conjectured to be bounded by one,
if $\pi$ is a supercuspidal representation (see \cite{hakim3} for details).
Such a pair is often called a supercuspidal Gelfand pair. It is known that
if ``almost all'' double cosets $HgH$ are fixed by an involution of $G$,
then $(G,H)$ is a supercuspidal Gelfand pair \cite{hakim5}.

It is natural to ask whether there is a symmetric space for which the 
multiplicity one property fails even for supercuspidal representations. There 
are such spaces, and an example is $(G={\rm SL}_n(E),H={\rm SL}_n(F))$, $n$ 
being an even
integer \cite{anand2}. 
This can be deduced from the formula for the dimension of 
${\rm Hom}_H(\pi,1)$ proved in \cite{anand2}. 
This is stated as Theorem 4.3 in this paper, and we also sketch a proof of it.
Also worth pointing out is 
the fact that the multiplicity formula in this context resembles closely to
the Labesse-Langlands multiplicity formula for the multiplicity of a 
representation in the cuspidal spectrum of ${\rm SL}_2({\Bbb A}_F)$ 
(for a number
field $F$) \cite{labesse,shelstad}.

Interestingly, when $n$ is odd, $({\rm SL}_n(E),{\rm SL}_n(F))$ is a 
Gelfand pair. We record this in the following theorem.
\begin{theorem}\label{gelfandpair}
Let $\pi$ be an irreducible admissible representation of ${\rm SL}_n(E)$. 
Then if $n$ is odd, $\dim_{\Bbb C}{\rm Hom}_{{\rm SL}_n(F)}(\pi,1)\leq 1$.
\end{theorem}

The plan of the paper is as follows. We recall some of the main results
in the case $({\rm GL}_n(E),{\rm GL}_n(F))$ in Section 2. The case 
$({\rm GL}_n(E),{\rm U}(n))$
is dealt with in the next section. Final two sections deal with the symmetric
pair $({\rm SL}_n(E),{\rm SL}_n(F))$, and the proof of Theorem 1.1. 
This is joint work with Dipendra Prasad.

It is a pleasure to thank R. Tandon (and the organisers) for inviting me to 
give a talk at the Hyderabad conference. I would also like to thank him for 
constant encouragement over the years. I would like to thank Dipendra Prasad 
for many discussions, suggestions, and encouraging words. Thanks are
also due to David Manderscheid, A. Raghuram, and
C.S. Rajan for questions and comments on this material. I would also
like to thank the referee for several helpful suggestions.

\section{The symmetric space $({\rm R}_{E/F}{\rm GL}(n),{\rm GL}(n))$}

Let $F$ be a $p$-adic field and let $E$ be a degree two extension of $F$.
Let $\sigma$ denote the nontrivial element of ${\rm Gal}(E/F)$. 
We recall some of the main theorems and conjectures about representations of
$G={\rm GL}_n(E)$ distinguished with respect to $H={\rm GL}_n(F)$. 

As already pointed out in the introduction $(G,H)$ is a Gelfand pair. The key
here is that there is an involution of $G$ that fixes all the double cosets
$HgH$. Specifically one knows the following lemma 
(Proposition 10, \cite{flicker1}).
\begin{prop}\label{involution}
The involution $g \mapsto g^{-\sigma}$ fixes the double cosets of $H$ 
in $G$.
\end{prop}

\begin{proof}[Proof (sketch)]
Consider the map from $G/H$ to the set 
$S=\{g \in G \mid gg^{\sigma}=
1\}$ given by $$gH \mapsto gg^{-\sigma}.$$ This map is a bijection
(it is clearly well defined and injective, and surjectivity follows 
by Hilbert 90). Next step is to show that two elements of $S$ which are 
conjugate in $G$ are conjugate in $H$. Now consider the elements
$gg^{-\sigma}$ and $g^{-\sigma}g$. Since these two elements are conjugate
in $G$, they should be conjugate in $H$. The proof follows from this. 
\end{proof}

Thus any distribution on $G$ that is $H$-bi-invariant is invariant under the
above involution. Also if an irreducible admissible representation $\pi$
of $G$ admits an $H$-invariant linear form, then so does the 
contragredient $\pi^\vee$ of $\pi$. Therefore the multiplicity one 
property follows from the well-known lemma due to Gelfand.

Proposition \ref{involution} also proves the following result 
(see Proposition 12, \cite{flicker1}).

\begin{theorem}\label{dist}
If $\pi$ is an irreducible admissible representation of $G$ which
is distinguished by $H$, then $\pi^\vee \cong \pi^\sigma$.
\end{theorem}

Also note that an obvious necessary condition for $\pi$ being distinguished
with respect to $H$ is that the central character $\omega_\pi$ of $\pi$
restricts trivially to $F^*$. Are these conditions sufficient too? A precise
conjectural answer is given below (here and elsewhere $\omega_{_{E/F}}$ 
signifies the quadratic character associated to the extension $E/F$).

\begin{conj}[Jacquet]\label{conj-jacquet}
Let $\pi$ be an irreducible admissible 
representation of ${\rm GL}_n(E)$ such that $\omega_\pi|_{_{F^*}}=1$ and
$\pi^\vee \cong \pi^\sigma$. Then $\pi$
is $H$-distinguished if $n$ is odd. If $n$ is even, $\pi$ is either
distinguished or $\omega_{_{E/F}}$-distinguished with respect to 
$H={\rm GL}_n(F)$.
\end{conj}

When $E/F$ is unramified, and $\pi$ a supercuspidal, this is proved by
D. Prasad \cite{dprasad3} where the proof eventually boils down to a similar
result about stable representations of finite groups of Lie type.
More recently the conjecture is settled by Anthony Kable \cite{kable}
whenever $\pi$ is a square integrable representation. The proof uses the theory
of the twisted tensor $L$-function (aka the Asai $L$-function) 
\cite{flicker2}, and has two main parts. First he proves the following
identity which relates the twisted $L$-function with the Rankin-Selberg
$L$-function.
\begin{align}\label{fact}
L(s,\pi \times \pi^\sigma)=L(s,As(\pi))L(s,As(\pi)\otimes\omega_{_{E/F}})
\end{align}
Then it is proved that if $L(s,As(\pi))$ has a pole at $s=0$, then $\pi$
is $H$-distinguished. Thus the Jacquet's conjecture follows
as it is well known that $L(s, \pi \times \pi^\vee)$ has a pole at 
$s=0$ (necessarily simple when $\pi$ is in the discrete series).

In fact the Asai $L$-function $L(s,As(\pi))$ having a pole at $s=0$
characterises $\pi$ being distinguished for a discrete series $\pi$.
The other direction is proved in \cite{anand3}. Thus we have the following
proposition.
\begin{prop}\label{dist-pole}
Let $\pi$ be a square integrable representation of ${\rm GL}_n(E)$. Then $\pi$
is distinguished with respect to ${\rm GL}_n(F)$ if and only if
$L(s,As(\pi))$ has a pole at $s=0$.
\end{prop}

Combining with the identity (\ref{fact}), this proves the following.
\begin{cor}\label{dist-cor}
Let $\pi$ be a square integrable representation of ${\rm GL}_n(E)$. Then
$\pi$ cannot be both distinguished and $\omega_{_{E/F}}$-distinguished
with respect to ${\rm GL}_n(F)$.
\end{cor}

Note that the central character considerations make this corollary obvious 
when $n$ is odd, and hence the nontrivial part is the case of an even $n$.

As already mentioned in the introduction, conjecturally, 
$H$-distinguished representations
lie in the image of the base change map from ${\rm U}(n,E/F)$ to 
${\rm GL}_n(E)$.
There are two base change maps from ${\rm U}(n)$ to ${\rm GL}_n(E)$, 
stable and unstable.
The precise conjecture is as follows (see \cite{flicker1}).

\begin{conj}[Flicker-Rallis]\label{conj-flicker}
Let $\pi$ be an irreducible admissible representation
of ${\rm GL}_n(E)$. Then $\pi$ is $H$-distinguished if and only if it is an
unstable (resp. stable) base change lift from ${\rm U}(n)$ if $n$ 
is even (resp. odd).
\end{conj}

Note that Conjecture \ref{conj-jacquet} follows from Conjecture
\ref{conj-flicker}.
When $\pi$ is in the discrete series, the images of the two base change maps
are believed to form two disjoint sets. Thus corollary \ref{dist-cor} is in 
agreement with this conjecture.

Dually we have (see \cite{jacquet1,jacquet2}):

\begin{conj}[Jacquet-Ye]\label{conj-ye}
The representation $\pi$ is distinguished with respect
to ${\rm U}(n)$ if and only if it is a base change lift from ${\rm GL}_n(F)$.
\end{conj}

\begin{remark}
A finite field analogue of the above two conjectures can be found
in \cite{gow}. Moreover it is proved there that both the pairs have the 
multiplicity one property (over a finite field).
\end{remark}

Conjecture \ref{conj-flicker} is known when $n=2$. 
In \cite{flicker1}, Flicker deduces it from
a similar global result. A local proof is provided in \cite{anand1}. There,
closely following a method due to Hiroshi Saito \cite{saito}, images of 
the base change
lifts are characterised in terms of the local factors of the representation,
which characterise distinguishedness by a result of Jeff Hakim \cite{hakim1}.
Recently, in a joint work with C.S. Rajan, we have been able to
prove Conjecture \ref{conj-flicker}
for $n=3$, when $\pi$ is in the discrete series. This is achieved
by closely analysing properties of the Asai $L$-function \cite{anand4}.
Conjecture \ref{conj-ye} (globally) is known for $n=2,3$ from 
the work of Jacquet and Ye.

\section{The symmetric space $({\rm R}_{E/F}{\rm U}(n),{\rm U}(n))$}

In this case $G={\rm GL}_n(E)$, and $H={\rm U}(n)$. The pair $(G,H)$ 
is conjectured to be
a supercuspidal Gelfand pair, and the conjecture has been proved in several 
cases (see \cite{hakim3}). In this section we only recall how this can be seen
when $n=2$.

Consider the group ${\rm GL}_2^+(F)=\{g \in {\rm GL}_2(F) 
\mid \det g \in N_{E/F}(E^*)\}$.
Then we know that $E^*{\rm GL}_2^+(F)=E^*{\rm U}(2)$. Using this identity one 
proves the following proposition.
\begin{prop}
Let $\pi$ be an irreducible
admissible representation of ${\rm GL}_2(E)$ with central character 
$\omega_\pi=
\mu\circ N_{E/F}$ for a character $\mu$ of $F^*$. Then $\pi$ is 
${\rm U}(2)$-distinguished if and only if it is either $\mu$-distinguished or
$\mu\omega_{_{E/F}}$-distinguished with respect to ${\rm GL}_2(F)$.
\end{prop} 
Moreover a $\mu$-distinguished (or $\mu\omega_{_{E/F}}$-distinguished)
functional is a ${\rm U}(2)$-distinguished functional. Conversely a 
${\rm U}(2)$-distinguished functional uniquely determines a nonzero functional
on which ${\rm GL}_2(F)$ acts by $\mu$ (or $\mu\omega_{_{E/F}}$).
Thus for a supercuspidal representation, the space of $H$-invariant linear 
forms has dimension less than or equal to one by Corollary \ref{dist-cor}.

On the other hand there are principal series representations of ${\rm GL}_2(E)$
which may admit a two dimensional space of $H$-invariant linear forms.
For instance, take $\pi=Ps(\chi,\chi^{-1})$ where $\chi$ is a character
of $E^*$ such that $\chi=\chi^\sigma$. Then $\pi$ is both distinguished
and $\omega_{_{E/F}}$-distinguished with respect to ${\rm GL}_2(F)$ 
\cite{hakim1}, and the corresponding functionals 
are ${\rm U}(2)$-distinguished.

\section{The symmetric space $({\rm R}_{E/F}{\rm SL}(n),{\rm SL}(n))$}

Representation theory of ${\rm SL}_n(E)$ can be understood in terms of the 
representation theory of ${\rm GL}_n(E)$ \cite{gelbart}. An 
irreducible admissible
representation $\pi$ of ${\rm SL}_n(E)$ comes in the restriction of 
an irreducible
admissible representation $\pi^\prime$ of ${\rm GL}_n(E)$. It is a theorem
of M. Tadic that the restriction to ${\rm SL}_n(E)$ of an irreducible 
admissible 
representation of ${\rm GL}_n(E)$ is multiplicity free \cite{tadic}. 
If two representations of ${\rm GL}_n(E)$ restrict to the same 
representation of ${\rm SL}_n(E)$,
then they are equivalent up to a twist by a character of $E^*$.
Constituents in the direct sum decomposition of a representation of 
${\rm GL}_n(E)$ form an $L$-packet of ${\rm SL}_n(E)$.

Let $\pi$ be an irreducible admissible representation of ${\rm SL}_n(E)$ which
is distinguished with respect to ${\rm SL}_n(F)$. Let $\pi^\prime$ be a 
representation of ${\rm GL}_n(E)$ such that $\pi$ occurs in the restriction of
$\pi^\prime$. Thus $\pi^\prime$ is distinguished with respect to 
${\rm SL}_n(F)$.
Consequently its central character $\omega_{\pi^\prime}$ is trivial
on the $n^{{\rm th}}$ roots of unity in $F^*$. It follows that $\omega_{\pi^
\prime}|_{_{F^*}}=\eta^n$ for a character $\eta$ of $F^*$. Thus we can assume,
after twisting by a character of $E^*$ if necessary, that $\omega_{\pi^\prime}$
restricts trivially to $F^*$. Now the space 
${\rm Hom}_{{\rm SL}_n(F)}(\pi^\prime,1)$
has the structure of a ${\rm GL}_n(F)$-module (with the obvious action) on 
which $F^*{\rm SL}_n(F)$ acts trivially. Thus it is a direct sum of characters 
of $F^*$ as a ${\rm GL}_n(F)$ module. Equivalently $\pi^\prime$ is 
$\chi$-distinguished with respect to ${\rm GL}_n(F)$ for a character 
$\chi$ of $F^*$.
Hence there is no loss of generality in assuming that $\pi^\prime$
is distinguished with respect to ${\rm GL}_n(F)$.
Also $\dim_{\Bbb C}{\rm Hom}_{{\rm SL}_n(F)}(\pi^\prime,1)$ equals the number
of characters $\chi$ of $F^*$ for which $\pi^\prime$ is $\chi$-distinguished
with respect to ${\rm SL}_n(F)$. This is so since 
$\dim_{\Bbb C}{\rm Hom}_{{\rm GL}_n(F)}
(\pi^\prime,\chi)\leq 1$.

Consider the subgroup of ${\rm GL}_n(E)$ defined by 
$${\rm GL}_n(E)^+=\{g \in {\rm GL}_n(E) \mid \det g \in F^*E^{*n}\}.$$ 
Also let us fix
a nontrivial additive character $\psi$ of $E$ which has trivial
restriction to $F$. Assume that $\pi^\prime$ is tempered (this is
because we need to use Proposition \ref{prop-pams}).

In the restriction of $\pi^\prime$ to ${\rm GL}_n(E)^+$, exactly 
one representation
is $\psi$-generic, say $\pi^+$. We contend that the ${\rm SL}_n(F)$-invariant
linear forms on $\pi^\prime$ are nontrivial only on the space of $\pi^+$
(among all the irreducible constituents of $\pi^\prime|_{_{{\rm GL}_n(E)^+}}$).
For the above claim we require the following proposition 
(Corollary 1.2, \cite{anand3}).
\begin{prop}\label{prop-pams}
Let $\pi^\prime$ be a tempered representation of ${\rm GL}_n(E)$ that is 
${\rm GL}_n(F)$-distinguished. Let $\psi$ be a nontrivial additive 
character of $E$
that has trivial restriction to $F$. Then the (nontrivial) 
${\rm GL}_n(F)$-invariant linear form on $\pi^\prime$ can be realized 
on the the 
$\psi$-Whittaker model of $\pi^\prime$
by $\ell(W)=\int_{N_n(F)\backslash P_n(F)}W(p)dp$. Here $P_n(F)$ is the
mirabolic subgroup of ${\rm GL}_n(F)$, and $N_n(F)$ is the unipotent radical
of the Borel subgroup of ${\rm GL}_n(F)$.
\end{prop}

Further we have,
\begin{prop}\label{prop-mrl}
All the constituents of the restriction of a representation of 
${\rm GL}_n(E)^+$
to ${\rm SL}_n(E)$ admit the same number of linearly independent 
${\rm SL}_n(F)$-invariant functionals.
\end{prop}

\begin{proof}
Indeed the constituents are conjugates of one another under the inner 
conjugation action of ${\rm GL}_n(F)$ on ${\rm SL}_n(F)$ 
(as ${\rm GL}_n(F){\rm SL}_n(E)E^*=
{\rm GL}_n(E)^+$). 
\end{proof}

Thus $\dim_{\Bbb C}{\rm Hom}_{{\rm SL}_n(F)}(\pi,1)$ is zero if $\pi$ does not
occur in the restriction of $\pi^+$ to ${\rm SL}_n(E)$. Moreover this dimension
is the same nonzero number for all $\pi$ that appear in the restriction
of $\pi^+$ to ${\rm SL}_n(E)$. From the preceding discussion 
it is clear that
this dimension is $$q(\pi)=q(\pi^\prime)=\frac{|X_{\pi^\prime}|}
{|Z_{\pi^\prime}|/
|Y_{\pi^\prime}|}$$
where 
$$X_{\pi^\prime}=\{\chi\in \hat{F^*}\mid \pi^\prime ~{\rm is}~ \chi-{\rm 
distinguished}\},$$
$$Y_{\pi^\prime}=\{\mu \in \hat{E^*}\mid \pi^\prime \otimes \mu 
\cong \pi^\prime, \mu|_{_{F^*}}=1\},$$
$$Z_{\pi^\prime}=\{\mu \in \hat{E^*}\mid \pi^\prime \otimes \mu \cong 
\pi^\prime\},$$
since $|Z_{\pi^\prime}|$ is the cardinality of the $L$-packet of $\pi$,
and $|Y_{\pi^\prime}|$ is the number of constituents in the restriction
of $\pi^\prime$ to ${\rm GL}_n(E)^+$ (by the result of Tadic 
quoted in the beginning
of this section). Note that the cardinalities of the sets above
do not depend on the choice of $\pi^\prime$.

Notice that $q(\pi)$ is a non-negative integer. This is so since the group
$Z_{\pi^\prime}/Y_{\pi^\prime}$ acts freely on $X_{\pi^\prime}$ (assume
$X_{\pi^\prime}$ is non-empty, $q(\pi)=0$ otherwise), and hence
$q(\pi)$ is the number of orbits.

All this information can be clubbed together as follows. Define a pairing
between $Z_{\pi^\prime}$ and the $L$-packet of $\pi$ by
$$<\mu,\pi>=\mu(a)$$ where $\pi$ is $\psi_a$-generic. ($\psi_a$ is the additive
character of $E$ given by $\psi_a(x)=\psi(ax)$). Then $$\frac{1}
{|Y_{\pi^\prime}|}\displaystyle{\sum_{\mu\in Y_{\pi^\prime}}}<\mu,\pi>$$
is one if $\pi$ comes in the restriction of $\pi^+$, and zero otherwise.
Thus we have (see Theorem 1.4, \cite{anand2}),
\begin{theorem}\label{thm-mrl}
Let $\pi$ be an irreducible admissible representation of 
${\rm SL}_n(E)$ that comes
in the restriction of a tempered representation of ${\rm GL}_n(E)$. Then,
$$\dim_{\Bbb C}{\rm Hom}_{{\rm SL}_n(F)}(\pi,1)=\frac{q(\pi)}{|Y_\pi|}
\displaystyle{\sum_{\mu\in Y_{\pi}}}<\mu,\pi>.$$
\end{theorem}

\begin{remark}
There is another way of looking at the number $q(\pi)$. Define
two equivalence relations on the set of all twists of $\pi^\prime$ that are
distinguished with respect to ${\rm GL}_n(F)$ as follows:
$$\pi_1^\prime \sim_w \pi_2^\prime \Longleftrightarrow \pi_1^\prime \cong
\pi_2^\prime \otimes \mu, \mu \in \hat{E^*}$$
$$\pi_1^\prime \sim_s \pi_2^\prime \Longleftrightarrow \pi_1^\prime \cong
\pi_2^\prime \otimes \mu, \mu \in \hat{E^*}, \mu|_{_{F^*}}=1.$$
Then $q(\pi)$ is the number of strong equivalence classes in a weak equivalence
class of $\pi^\prime$.
\end{remark}

\section{Proof of Theorem \ref{gelfandpair}}

Let $\pi$ be an irreducible admissible representation of ${\rm SL}_n(E)$ which
is distinguished with respect to ${\rm SL}_n(F)$. Let $\pi^\prime$ be a 
representation of ${\rm GL}_n(E)$ such that $\pi$ comes in the restriction of
$\pi^\prime$ to ${\rm SL}_n(F)$. 

We need to introduce the group
$$Y^\prime_{\pi^\prime}=\{\mu \in \hat{E^*} \mid \pi^\prime \otimes \mu \cong
\pi^\prime, \mu|_{_{E^1}}=1\}$$
where $E^1$ denotes the norm one elements of $E^*$.
We claim that the map
$$\chi \mapsto \chi \circ N_{E/F}$$
establishes an injection from $X_{\pi^\prime}$ to $Y^\prime_{\pi^\prime}$. 

This is of course vacuously true if $X_{\pi^\prime}$ is empty. Otherwise, as
before, we can assume that $\pi^\prime$ is distinguished with respect to
${\rm GL}_n(F)$. 
Then this is a well defined map follows from an application of Theorem 
\ref{dist}. 
Injectivity is a consequence of the fact that $\pi^\prime$ cannot be both
distinguished and $\omega_{_{E/F}}$-distinguished when $n$ is an odd integer.
(This can be seen by considering central characters. This fact is not true
in general when $n$ is even. See also Corollary \ref{dist-cor}).
If we assume the truth of Jacquet's conjecture, this map is in fact 
a bijection, but we do not need that.

Thus $$q(\pi)\leq \frac{|Y_{\pi^\prime}||Y^\prime_{\pi^\prime}|}
{|Z_{\pi^\prime}|}.$$ If $n$ is an odd integer, note that $Y_{\pi^\prime}
\bigcap Y^\prime_{\pi^\prime}=\{1\}$. Therefore the quantity on the right
side of the above inequality is 
$\frac{|Y_{\pi^\prime}Y^\prime_{\pi^\prime}|}
{|Z_{\pi^\prime}|}$. This forces $q(\pi)$ to be zero or one, since we 
know that it is a non-negative integer.

Now let $$\pi^\prime|_{_{{\rm GL}_n(E)^+}}=\oplus_i \pi_i^+$$
be the direct sum decomposition of $\pi^\prime$ 
restricted to ${\rm GL}_n(E)^+$, where 
$\pi_i^+$ are inequivalent. Let $a_i$ denote the common dimension
(by Proposition \ref{prop-mrl}) of the space of 
${\rm SL}_n(F)$-invariant forms 
on constituents
of $\pi_i^+|_{_{{\rm SL}_n(E)}}$. We conclude that $$\sum_i a_i=q(\pi).$$  
Thus $a_i$ cannot be more than one, proving the theorem.

When $n$ is an even integer, $Y_{\pi^\prime}\bigcap Y^\prime_{\pi^\prime}
\neq \{1\}$ in general, and this results in higher multiplicities. 
We restrict ourselves to just one example. 

To this end, consider a quadratic extension $K$ of $F$ different from $E$,
and let $\eta$ be a character of $K^*$ with trivial restriction to $F^*$.
Let $L$ denote the compositum of $E$ and $K$. Also assume that $K/F$
is such that we can choose $\eta$ so that $\eta^8 \neq 1$.
Let $\pi_0$ be the representation of ${\rm GL}_2(F)$
obtained by automorphically inducing $\eta$, and let $\pi^\prime$ be the
base change lift of $\pi_0$ to ${\rm GL}_2(E)$.
 
Our assumption on $\eta$ guarantees that $\pi^\prime$ is a supercuspidal 
representation, and that $|Z_{\pi^\prime}|=2$ (see \cite{labesse,shelstad}). 
This forces $|Y^\prime_{\pi^\prime}|=2$ (as either $\mu$ or $\mu\circ N_{E/F}$
has to be in $Y^\prime_{\pi^\prime}$ if $\mu \in Z_{\pi^\prime}$).
Since $n=2$, $Y_{\pi^\prime}=Y^\prime_{\pi^\prime}$. Thus in this case
$q(\pi)=2$.

\end{document}